\documentclass[12pt]{article}
\usepackage[dvips]{epsfig}
\usepackage{amsmath}
\usepackage{amsfonts}
\usepackage{amssymb}
\usepackage{graphicx}

\numberwithin{equation}{section}

\begin{document}

\title{On the nodal line of the second eigenfunction of
the Laplacian over some concave domains in $\mathbb{R}^2$}
\author{Donghui Yang \\
School of Mathematics and Statistics, Central China Normal University,\\ 430079, Wuhan, People's Republic of China
\\ \textsf{donghyang@mail.ccnu.edu.cn}\\}
\date{}
\maketitle

\begin{abstract}
In this paper we will prove the nodal line $N$ of the second
eigenfunction of the Laplacian over some simply connected concave
domain $\Omega$ in $\mathbb{R}^2$ must intersect the boundary
$\partial\Omega$ at exactly two points.

\end{abstract}

\thispagestyle{empty}


\baselineskip=17.5pt \parskip=3pt

\section{ Introduction}

\hskip600pt \footnote{The authors were supported by NSFC 10901069.}

An eigenfunction $\varphi_i$ is meant to be a solution of
Dirichlet's problem:
\begin{eqnarray}\label{1.1}
\begin{cases}  \Delta \varphi_i+\lambda_i \varphi_i=0
\quad &\mbox{in}\ \Omega \cr \varphi_i=0\quad &\mbox{on}\
\partial\Omega,
\end{cases}
\end{eqnarray}
where $\displaystyle\Delta=\frac{\partial^2}{\partial
x_1^2}+\frac{\partial^2}{\partial x_2^2}$ is the Laplacian, $\Omega$
is a bounded smooth domain in $\mathbb{R}^2$, $\lambda_i$ is the
$i$th eigenvalue with $\lambda_1<\lambda_2\leq\lambda_3\leq\cdots$,
and $\varphi_i$ is the $i$th eigenfunction ($i=1, 2, \cdots$). It is
well know that the first eigenfunction is positive in $\Omega$, and
all higher eigenfunctions must change sign. The nodal set of an
eigenfunction $\varphi_i$ is defined to be the closure of
$\{x\in\Omega;\ \varphi_i(x)=0\}$. The Courant nodal domain theorem
\cite{C1} tells us that the nodal set of an $i$th eigenfunction $\varphi_i$
divides the domain $\Omega$ into at most $i$ subregions. especially, $\varphi_2$ divides the
domain $\Omega$ into at exactly 2 domains.

In 1967 Payne \cite{Pa2} conjectured that $\varphi_2$ cannot have a
closed nodal line in $\Omega$ and in 1982 Yau \cite{Yau} asked the same
question for convex domains in $\mathbb{R}^2$. Payne \cite{Pa1} proved that the nodal line
touches the boundary of a convex set which is symmetric under a reflection.
C.-S. Lin \cite{L1} proved the conjecture
provided the domain $\Omega\subset\mathbb{R}^2$ is smooth, convex, and invariant under
a rotation with angle $\displaystyle\frac{2\pi p}{q}$, where $p$ and $q$ are positive
integers. D. Jerison \cite{J1} proved the conjecture for long thin convex sets.
 Melas \cite{Me} have
settled the convex case for $C^\infty$ boundary and this was
extended to general boundary by Alessandrini \cite{Al}.
M. Hoffmann-Ostenhof, T. M. Hoffmann-Ostenhof and N. Nadirashvili \cite{H1}
construct a nonconvex, not simply connected domain for which the second
eigenfunction has a closed nodal line. Also for convex $D$ D. Jerison \cite{J2}
and D. Grieser and D. Jerison \cite{G1} obtained interesting results on the location
of the first nodal line.

In this paper we obtain that the nodal line of the second
eigenfunction $\varphi_2$ over some simply connected concave domains
$\Omega$ intersect the boundary $\partial\Omega$ at exactly two
points. Which is the special case of the following theorem:

\noindent{\bf Theorem }\quad The nodal line of a second eigenfunction of Laplacian divides the domain $\Omega$
by intersecting its boundary at exactly two points if  the domain $\rho(\Omega)$ is strictly convex in $\theta$ and
symmetric with respect to the $r$-axis.

\section {main results and their proofs}

Let $\Omega,\widetilde\Omega\subset\mathbb{R}^2$ be two smooth domains (0 is not in the closure of
$\widetilde\Omega$) and
$\rho: \Omega\rightarrow \widetilde\Omega,\ (x, y)\mapsto (r, \theta)$ be a diffeomorphism
defined by $(x, y)=\rho^{-1}(r, \theta)$ with $x=r\cos\theta, y=r\sin\theta$. Then the equation
(\ref{1.1}) becomes a new equation:
\begin{eqnarray}\label{1.2}
\begin{cases}\displaystyle \frac{\partial^2 \widetilde\varphi_i}{\partial r^2}
+\frac{1}{r^2}\frac{\partial^2 \widetilde\varphi_i}{\partial \theta^2}
+\frac{1}{r}\frac{\partial \widetilde\varphi_i}{\partial r}+\lambda_i \widetilde\varphi_i=0 & \mbox{in}\ \widetilde{\Omega} \cr
\widetilde\varphi_i=0 &\mbox{on}\ \partial\widetilde{\Omega},
\end{cases}
\end{eqnarray}
where $\widetilde\varphi_i(r, \theta)=\varphi_i\circ\rho^{-1}(r, \theta)$.

Now we display the relations between the equations (\ref{1.1}) and (\ref{1.2}): If $\varphi_i$ is
a solution of equation  (\ref{1.1}), then $\varphi_i$ is a smooth function and
 $\widetilde\varphi_i=\varphi_i\circ\rho^{-1}$
is also a smooth solution of (\ref{1.2}); conversely, if $\widetilde\varphi_i$ is a solution of (\ref{1.2}), then
$\widetilde\varphi_i\in H_0^1(\widetilde\Omega)$, and $\widetilde\varphi_i\in C^\infty(\widetilde\Omega)$
by the infinite differentiability up to the boundary theorem (See \cite{E1} pp. 324-326), hence
$\varphi_i=\widetilde\varphi_i\circ\rho$ is a solution of equation (\ref{1.1}). And hence
$\lambda_2(\Omega)=\lambda_2(\rho(\Omega))$.

Following from Courant nodal domain theorem
\cite{C1} we know that the nodal set of an $i$th eigenfunction
$\widetilde\varphi_i$ divides
the domain $\widetilde\Omega$ into at most $i$ subregions, especially,
$\widetilde\varphi_2$ divides the domain
$\widetilde\Omega$ into at exactly $2$ subregions.

Throughout the paper we denote $\varphi_2$ an second eigenfunction of (\ref{1.1}) and $N=\overline{\{x\in\Omega;\
\varphi_2(x)=0\}}$ is the nodal line of $\varphi_2$; $\widetilde\varphi_2=\varphi_2\circ\rho^{-1}$ and
$\widetilde N=\overline{\{x\in\widetilde\Omega;\ \widetilde\varphi_2(x)=0\}}$.

The following Lemma 1 is proved in \cite{L1}:

\vskip 2mm

\noindent{\bf Lemma 1}\quad Suppose $ P\in \partial \Omega$. Then $\displaystyle
\frac{\partial\varphi_2}{\partial \nu}( P)=0$ if and only if
$P \in N$, where
$\displaystyle\frac{\partial\varphi_2}{\partial \nu}$ is the outnormal
derivative of $\varphi_2$ on the boundary.

\vskip 2mm

\noindent{\bf Proof}\quad Now, we prove Lemma 1 again as some different way.

Let $P\in\partial\Omega$. Since $\varphi_2$ is a smooth solution of equation (\ref{1.1}) and
$\Omega$ is a smoothly bounded domain, there exists an open set $W\subset\mathbb{R}^2$
(See \cite{E1}, pp. 254-256) such that
$\varphi^*_2$ is an extension of $\varphi_2$ in $W$, and
\begin{eqnarray}\label{1.3}
P\in W,\quad
\varphi^*_2\in C^1(W)\quad \mbox{and}\quad
\varphi^*_2|_{\Omega\cap W}=\varphi_2, \quad
\varphi^*_2|_{\partial\Omega\cap W}=0.
\end{eqnarray}

Let $P=(x_0,y_0) \in \partial\Omega\cap N$. If $\displaystyle\frac{\partial\varphi_2}{\partial\nu}(P)\not=0$,
then $\displaystyle\pm|\nabla\varphi^*_2|(P)=\nabla\varphi^*_2\cdot\nu(P)
=\frac{\partial\varphi^*_2}{\partial\nu}(P)=\frac{\partial\varphi_2}{\partial\nu}(P)\not=0$
(Please see Appendix for the first equality). Without loss of generality, we assume
$\displaystyle \frac{\partial\varphi^*_2}{\partial y}(P)\not=0$.
By Implicit Function Theorem, there exists a unique function
$g: (x_0-\epsilon, x_0+\epsilon)\rightarrow \mathbb{R}$ ($\epsilon$ is small enough)
such that $g\in C^1(x_0-\epsilon, x_0+\epsilon)$, $g(x_0)=y_0$ and for any $x\in (x_0-\epsilon, x_0+\epsilon)$ we have
$\varphi^*_2(x,g(x))=\varphi^*_2(x_0,y_0)=0$. Which implies that $(x, g(x))\in\partial\Omega$ by (\ref{1.3}),
hence $\varphi^*_2(x,y)\not=0$ for $(x,y)\in\Omega\cap W$. Which contradict to $P\in N$. i.e.
$\displaystyle\frac{\partial\varphi_2}{\partial\nu}(P)=\pm|\nabla\varphi^*_2|(P)=0$.

Secondly, suppose $P\in N$, then, by the Hopf's Lemma, we have $\displaystyle
\frac{\partial\varphi_2}{\partial\nu}(P)\not=0$.\quad $\Box$

\vskip 2mm

By the same way, we obtain the following Lemma 2:

\vskip 2mm

\noindent{\bf Lemma 2}\quad Suppose $\widetilde P\in \partial\widetilde\Omega$. Then $\displaystyle
\frac{\partial\widetilde\varphi_2}{\partial\widetilde \nu}(\widetilde P)=0$ if and only if
$\widetilde P\in \widetilde N$, where
$\displaystyle\frac{\partial\widetilde\varphi_2}{\partial\widetilde\nu}$ is the outnormal
derivative of $\widetilde\varphi_2$ on the boundary.

\vskip 2mm

\noindent{\bf Proof}\quad Let $\widetilde P\in\partial\widetilde\Omega$. Since $\widetilde\varphi_2$ is a smooth solution of equation (\ref{1.2}) and
$\widetilde\Omega$ is a smoothly bounded domain, there exists an open set $\widetilde W\subset\mathbb{R}^2$ such that
$\widetilde\varphi^*_2$ is an extension of $\widetilde\varphi_2$ in $\widetilde W$, and
\begin{eqnarray}\label{1.3}
\widetilde P\in \widetilde W,\quad
\widetilde\varphi^*_2\in C^1(\widetilde W)\quad \mbox{and}\quad
\widetilde\varphi^*_2|_{\widetilde\Omega\cap \widetilde W}=\widetilde\varphi_2, \quad
\widetilde\varphi^*_2|_{\partial\widetilde\Omega\cap \widetilde W}=0.
\end{eqnarray}

Let $\widetilde P=(r_0, \theta_0) \in \partial\widetilde\Omega\cap \widetilde N$. If $\displaystyle\frac{\partial\widetilde\varphi_2}{\partial\widetilde\nu}(\widetilde P)\not=0$,
then $\displaystyle\pm|\nabla\widetilde\varphi^*_2|(\widetilde P)=\nabla\widetilde\varphi^*_2\cdot\widetilde\nu(\widetilde P)
=\frac{\partial\widetilde\varphi^*_2}{\partial\widetilde\nu}(\widetilde P)=\frac{\partial\widetilde\varphi_2}{\partial\widetilde\nu}(\widetilde P)\not=0$. Without loss of generality, we assume
$\displaystyle \frac{\partial\widetilde\varphi^*_2}{\partial \theta}(\widetilde P)\not=0$.
By Implicit Function Theorem, there exists a unique function
$g: (r_0-\epsilon, r_0+\epsilon)\rightarrow \mathbb{R}$ ($\epsilon$ is small enough)
such that $g\in C^1(r_0-\epsilon, r_0+\epsilon)$, $g(r_0)=\theta_0$ and for any $r\in (r_0-\epsilon,  r_0+\epsilon)$ we have
$\widetilde\varphi^*_2(r,g(r))=\widetilde\varphi^*_2(r_0,\theta_0)=0$. Which implies that $(r, g(r))\in\partial\widetilde\Omega$ by (\ref{1.3}),
hence $\widetilde\varphi^*_2(r,\theta)\not=0$ for $(r,\theta)\in\widetilde\Omega\cap \widetilde W$. Which contradict to $\widetilde P\in\widetilde N$. i.e.
$\displaystyle\frac{\partial\widetilde\varphi_2}{\partial\widetilde\nu}(\widetilde P)=\pm|\nabla\widetilde\varphi^*_2|(\widetilde P)=0$.

Now, we show that if $\displaystyle
\frac{\partial\widetilde\varphi_2}{\partial\widetilde \nu}(\widetilde P)=0$ then
$\widetilde P\in \widetilde N$.

Otherwise, there exist a neighborhood $U\subset\mathbb{R}^2$ of $\widetilde P$ such that
$N\cap U=\emptyset$. Then either $\widetilde\varphi_2(r,\theta)>0$ or $\widetilde\varphi_2(r,\theta)<0$
for all $(r,\theta)\in U\cap\widetilde\Omega$, without loss of generality, we assume that
$\widetilde\varphi_2(r,\theta)<0$ for all $(r,\theta)\in U\cap\widetilde\Omega$.
Since $\widetilde\varphi_2$ divides the domain
$\widetilde\Omega$ into at exactly $2$ subregions $\widetilde\Omega_1$ and $\widetilde\Omega_2$,
we take $U\cap\widetilde\Omega
\subset \widetilde\Omega_1$, then $\widetilde\varphi_2(r,\theta)<0$ for all
$(r,\theta)\in\widetilde\Omega_1$. Since $\displaystyle{\Big(}-\frac{\partial^2 }{\partial r^2}
-\frac{1}{r^2}\frac{\partial^2 }{\partial \theta^2}
-\frac{1}{r}\frac{\partial }{\partial r}{\Big)}(\widetilde\varphi_2)=\lambda_2 \widetilde\varphi_2\leq0$
in $\widetilde\Omega_1$, and $\widetilde\varphi_2(\widetilde P)=0>\widetilde\varphi_2(r,\theta)$ for all
$(r,\theta)\in\widetilde\Omega_1$, we obtain $\displaystyle
\frac{\partial\widetilde\varphi_2}{\partial\widetilde \nu}(\widetilde P)>0$ by Hopf's Lemma
(See \cite{E1} pp. 330-332). Which implies the contradiction. \quad $\Box$

\vskip 2mm

Now we consider the equation:
\begin{eqnarray}\label{1.4}
\begin{cases}\displaystyle \frac{\partial^2 \widetilde\varphi_2}{\partial r^2}
+\frac{1}{r^2}\frac{\partial^2 \widetilde\varphi_2}{\partial \theta^2}
+\frac{1}{r}\frac{\partial \widetilde\varphi_2}{\partial r}+\lambda_2 \widetilde\varphi_2=0 & \mbox{in}\ \widetilde{\Omega} \cr
\widetilde\varphi_2=0 &\mbox{on}\ \partial\widetilde{\Omega},
\end{cases}
\end{eqnarray}
where $\widetilde{\Omega}$ is a smoothly bounded domain, $\widetilde{\Omega}\subset\{(r, \theta)
\in\mathbb{R}^2;\ 0<r_0\leq r\leq R_0\}$ ($r_0, R_0$ are given constants),
 $\widetilde{\Omega}$ is strictly convex in $\theta$ and $\widetilde{\Omega}$ is symmetric with respect to $r$-axis. (We note that
the domain $\widetilde\Omega$ is strictly convex in $\theta$ if every line parallel to the $\theta$-axis which intersect $\widetilde\Omega$, cut
$\partial\Omega$ in at most two points.)

Then we obtain the following Lemma 3:

\vskip 2mm

\noindent{\bf Lemma 3 }\quad The nodal line $\widetilde N$ of a second eigenfunction $\widetilde\varphi_2$ divides the domain
$\widetilde\Omega$ by intersecting its boundary at exactly two points.

\vskip 2mm

\noindent{\bf Proof}\quad The proof of this lemma is similar to Theorem 2.2 in \cite{L1} and Theorem I in \cite{Pa2}.

If $\widetilde \varphi_2$ is odd in $\theta$ (i.e. $\widetilde \varphi_2(r,-\theta)=-\widetilde \varphi_2(r,\theta)$),
the nodal line $\widetilde N$ is just the $r$-axis. And Lemma 3 is obviously true. Suppose $\widetilde \varphi_2$ is
even in $r$ (i.e. $\widetilde \varphi_2(r,-\theta)=\widetilde \varphi_2(r,\theta)$). Assume that Lemma 3 is false.
Then for $\widetilde P\in\partial\widetilde\Omega$, $\displaystyle(\frac{\partial\widetilde\varphi}{\partial \theta})
(\widetilde P)\not=0$, except the tangent of $\partial\Omega$ at $\widetilde P$ is the $\theta$-direction. Without loss
of generality, we may assume that $\displaystyle(\frac{\partial\widetilde\varphi}{\partial \theta})
(\widetilde P)\geq 0$ for $\widetilde P\in\partial\widetilde\Omega\cap\{(r,\theta);\ \theta\leq0\}$. Set $(\widetilde
\Omega)^-=\mbox{Cl}(\widetilde\Omega)\cap\{(r,\theta);\ \theta\leq0\}$, here and throughout the paper
$\mbox{Cl}(\widetilde\Omega)$ represent the
closure of $\widetilde\Omega$. $\displaystyle \frac{\partial\widetilde\varphi}{\partial \theta}$ must change sign in
$(\widetilde\Omega)^-$. Otherwise $\widetilde \varphi_2\geq0$ in $(\widetilde\Omega)^-$. by evenness,
$\widetilde \varphi_2\geq0$ in $\mbox{Cl}(\widetilde\Omega)$, which leads to a contradiction. Hence the nodal line
$\displaystyle\{(r,\theta)\in(\widetilde\Omega)^-;\  \frac{\partial\widetilde\varphi}{\partial \theta}(r,\theta)=0\}$
encloses a subregion $(\widetilde\Omega)^-_*$ of
$\displaystyle\{(r,\theta)\in(\widetilde\Omega)^-;\  \frac{\partial\widetilde\varphi}{\partial \theta}(r,\theta)<0\}$.
Let $\widetilde\Omega_*=\{(r,\theta)\in\widetilde\Omega;\ \mbox{either}\ (r, \theta)\in(\widetilde\Omega)^-_*\ \mbox{or}\
(r,-\theta)\in(\widetilde\Omega)^-_*\}$. Then $\displaystyle\frac{\partial\widetilde\varphi_2}{\partial \theta}$ satisfies
\begin{eqnarray}\label{1.5}
&&{\Big(}\frac{\partial^2 }{\partial r^2}
+\frac{1}{r^2}\frac{\partial^2 }{\partial \theta^2}
+\frac{1}{r}\frac{\partial }{\partial r}+\lambda_2{\Big)}(\frac{\partial\widetilde\varphi_2}{\partial \theta})\cr
&&=\frac{\partial}{\partial \theta}{\Big[}{\Big(}\frac{\partial^2 }{\partial r^2}
+\frac{1}{r^2}\frac{\partial^2 }{\partial \theta^2}
+\frac{1}{r}\frac{\partial }{\partial r}+\lambda_2{\Big)}(\widetilde\varphi_2){\Big]}=0 \quad \mbox{in}\ \widetilde\Omega_*,
\end{eqnarray}
and
\begin{eqnarray}\label{1.6}
\displaystyle\frac{\partial\widetilde\varphi_2}{\partial \theta}=0\quad \mbox{on}\ \partial\widetilde\Omega_*.
\end{eqnarray}
By (\ref{1.5}) and (\ref{1.6}) we obtain that
$\displaystyle\widehat{\varphi}_2:=\frac{\partial\widetilde\varphi_2}{\partial \theta}\circ\rho$ satisfies
\begin{eqnarray*}
\begin{cases}  \Delta \widehat\varphi_2+\lambda_2 \widehat\varphi_2=0
\quad &\mbox{in}\ \Omega_*:=\rho^{-1}(\widetilde\Omega_*) \cr
\widehat\varphi_2=0\quad &\mbox{on}\ \partial\Omega_*,
\end{cases}
\end{eqnarray*}

Since $\displaystyle \frac{\partial\widetilde\varphi_2}{\partial \theta}(r,-\theta)
=-\frac{\partial\widetilde\varphi_2}{\partial \theta}(r,\theta)$,
$\displaystyle\frac{\partial\widetilde\varphi_2}{\partial \theta}$
must change sign, Then $\widehat \varphi_2$ is also change sign.
Hence $\lambda_2\geq \lambda_2(\Omega_*)$, where $\lambda_2(\Omega_*)$ is the second
eigenvalue of Laplacian in $\Omega_*$.
But $ \Omega_*\subset \Omega$ (since $\widetilde\Omega_*\subset\widetilde\Omega$), and by monotony principle
we get
\begin{eqnarray*}
\lambda_2\geq\lambda_2(\Omega_*)>\lambda_2,
\end{eqnarray*}
which is a contradiction.

In general a second eigenfunction $\widetilde\varphi_2$ can be written as $\widetilde\varphi_2=\phi_1+\phi_2$, where
$\displaystyle\phi_1(r,\theta)=\frac{\widetilde\varphi_2(r,\theta)-\widetilde\varphi_2(r,-\theta)}{2}$
is odd in $\theta$,
$\displaystyle\phi_2(r,\theta)=\frac{\widetilde\varphi_2(r,\theta)+\widetilde\varphi_2(r,-\theta)}{2}$
is even in $\theta$ and both are second eigenfunctions. By the above proof, we know that there exist two points
$\widetilde P=(r_0, \theta_0)$ and $\widetilde Q=(r_0, -\theta_0)\in\partial\widetilde\Omega$, where $\theta_0\not=0$
such that $\displaystyle(\frac{\partial\phi_2}{\partial \widetilde\nu})(\widetilde P)=
(\frac{\partial\phi_2}{\partial \widetilde\nu})(\widetilde Q)=0$. Since
$\displaystyle\frac{\partial\phi_1}{\partial \widetilde\nu}$ is odd on $\partial\widetilde\Omega$, we assume
$\displaystyle(\frac{\partial\phi_1}{\partial \widetilde\nu})(\widetilde P)>0$ and
$\displaystyle(\frac{\partial\phi_1}{\partial \widetilde\nu})(\widetilde Q)<0$. Hence
$\displaystyle(\frac{\partial\widetilde\varphi_2}{\partial \widetilde\nu})(\widetilde P)>0$ and
$\displaystyle(\frac{\partial\widetilde\varphi_2}{\partial \widetilde\nu})(\widetilde Q)<0$. It implies
that there exist two points $\widetilde P^*$ and $\widetilde Q^*$ on $\partial\widetilde\Omega$ such that
\begin{eqnarray*}
(\frac{\partial\phi_2}{\partial \widetilde\nu})(\widetilde P^*)=
(\frac{\partial\phi_2}{\partial \widetilde\nu})(\widetilde Q^*)=0,
\end{eqnarray*}
and the Lemma 3 for this case follows. \quad $\Box$

\vskip 2mm

Since $\rho: \Omega\rightarrow \widetilde\Omega,\ (x, y)\mapsto (r, \theta)$ is a diffeomorphism
defined by $(x, y)=\rho^{-1}(r, \theta)$ with $x=r\cos\theta, y=r\sin\theta$, we get the following theorem:

\vskip 2mm

\noindent{\bf Theorem }\quad The nodal line of a second eigenfunction of Laplacian divides the domain $\Omega$
by intersecting its boundary at exactly two points if  the domain $\rho(\Omega)$ is strictly convex in $\theta$ and
symmetric with respect to the $r$-axis.

\vskip 2mm

\noindent{\bf Remark }\quad At the end of this paper, we give a figure,
which displays that the diffeomorphism $\rho(\Omega)$ of
some concave domain $\Omega$ satisfies the conditions of Theorem.

In the following figure, $\Omega$ in figure (A) is the domain formed by the curves $C_1, C_2, C_3, C_4$, where
$\displaystyle C_1=\{(x,y)\in\mathbb{R}^2;\ \frac{x^2}{6^2}+\frac{y^2}{3^2}=1\ \mbox{and}\ x\geq 0\}$,
$\displaystyle C_2=\{(x,y)\in\mathbb{R}^2;\ \frac{x^2}{1^2}+\frac{y^2}{2^2}=1\ \mbox{and}\ x\geq 0\}$,
$\displaystyle C_3=\{(x,y)\in\mathbb{R}^2;\ x^2+(y-\frac{5}{2})^2=(\frac{1}{2})^2\ \mbox{and}\ x\leq 0\}$,
$\displaystyle C_4=\{(x,y)\in\mathbb{R}^2;\ x^2+(y+\frac{5}{2})^2=(\frac{1}{2})^2\ \mbox{and}\ x\leq 0\}$.
$\rho(\Omega)$ in figure (B) is the diffeomorphism of $\Omega$. Obviously, $\Omega$ is a concave domain and $\rho(\Omega)$ is a convex domain.

\begin{figure}[h]
\begin{center}
\resizebox{0.9\textwidth}{!}{%
\includegraphics{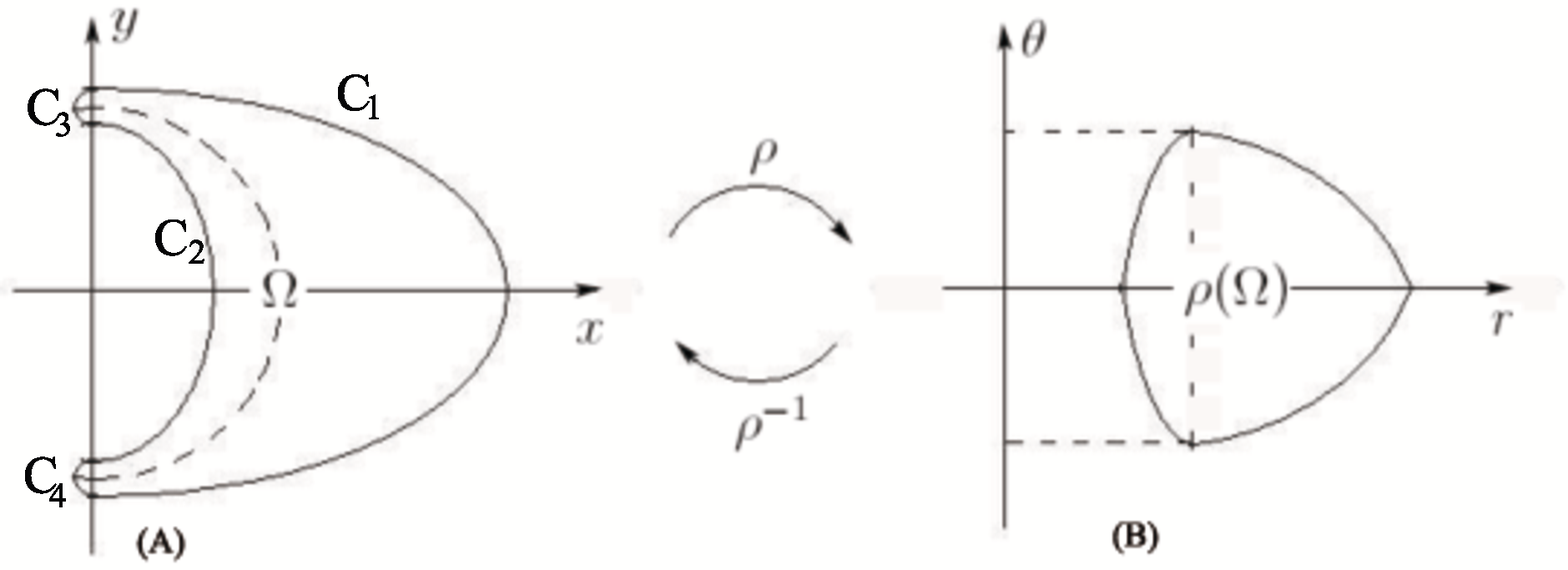}}
\end{center}
\end{figure}

\end{document}